\documentclass{article}

\usepackage{multirow}
\usepackage{graphicx}
\usepackage{url}
\newcommand{\figscale}{0.79}

\title{A$^2$ILU: Auto-accelerated ILU Preconditioner for Sparse Linear Systems\thanks{}}

\author{Yuichiro Miki\thanks{Simulation \& Analysis R\&D Center, Canon Inc., 30-2, Shimomaruko 3-chome, Ohta-ku, Tokyo 146-8501, Japan ({\tt miki.yuichiro@canon.co.jp}).} \and Teruyoshi Washizawa\thanks{Simulation \& Analysis R\&D Center, Canon Inc., 30-2, Shimomaruko 3-chome, Ohta-ku, Tokyo 146-8501, Japan ({\tt washizawa.teruyoshi@canon.co.jp}).}}

\date{}

\begin{document}

\maketitle

\begin{abstract}
The ILU-based preconditioning methods in previous work have their own parameters to improve their performances.
Although the parameters may degrade the performance, their determination is left to users.
Thus, these previous methods are not reliable in practical computer-aided engineering use.
This paper proposes a novel ILU-based preconditioner called the auto-accelerated ILU, or A$^2$ILU.
In order to improve the convergence, A$^2$ILU introduces acceleration parameters which modify the ILU factorized preconditioning matrix.
A$^2$ILU needs no more operations than the original ILU because the acceleration parameters are optimized automatically by A$^2$ILU itself.
Numerical tests reveal the performance of A$^2$ILU is superior to previous ILU-based methods with manually optimized parameters.
The numerical tests also demonstrate the ability to apply auto-acceleration to ILU-based methods to improve their performances and robustness of parameter sensitivities.
\end{abstract}



\pagestyle{myheadings}
\thispagestyle{plain}
\markboth{Y. MIKI AND T. WASHIZAWA}{A$^2$ILU: AUTO-ACCELERATED ILU PRECONDITIONER}

\section{Introduction}
As a means of fast solving ``sparse'' linear systems which can be the most time-consuming part in many physical simulations, a preconditioned iterative method is one of the most frequently used computational methods.
Among many existing preconditioning methods, ILU preconditioning is highly regarded for its generality in application because it can be applied to arbitrary matrices with nonzero entry structures.
In particular ILU with no fill-in, denoted by ILU(0) \cite{ILU}, which is the most basic form of ILU preconditioning, requires no other information except the equations themselves, i.e., it is a parameter-free method, making it extremely practical.
Because of these properties, it has achieved success in a wide variety of physical simulations.

Several ILU-based methods have been proposed to improve the computational performance of ILU(0) \cite{Templates, Benzi, NLAforHPC, Saad}.
These are classified into two types: one changes only the nonzero entry values in the preconditioning matrices obtained by ILU(0) and the other changes both the nonzero entry structure and their values.
The first type includes shifted ILU \cite{shiftedILU} and modified ILU (MILU) \cite{modifiedILU, Gustafsson}.
Shifted ILU performs ILU factorization on a matrix obtained by shifting the diagonal entries.
MILU subtracts the products of a relaxation parameter and fill-ins to be discarded from the diagonal entries to reduce the approximation error.
The second type includes ILU with $k$ extra diagonals \cite{Gustafsson}, fill-in level ILU \cite{fill-inLevelILU}, and ILUT \cite{ILUT}, as well as Crout ILU \cite{CroutILU}.
The ILU with $k$ extra diagonals chooses whether to allow or discard fill-ins, depending on the position in the discrete space.
The fill-in level ILU is the extended version of this method for a general sparse matrix.
ILUT and Crout ILU allow fill-ins only within certain ranges set up for the number of nonzero entries and those values.
The first type of method can be applied to other ILU-based methods.

Each of these ILU-based methods obtains the improvement of the performance by introducing unique parameters.
The performance is degraded, however, by setting inadequate parameter values.
Since none of these methods sets up the parameters automatically, this should be left to users.
A straightforward way to optimize the parameters is brute force searching over a set of candidates.
Here we denote the number of candidates as $N_S$.
This type of optimization needs to solve a given linear equation $N_S$ times.
This approach is therefore available only when we solve the equation with the same coefficient matrix for a different right-hand-side vector much more than $N_S$ times.

At workplaces where computer-aided engineering is used, a wide variety of physical simulations is used to describe various physical phenomena by different modeling methods.
Even for one physical simulation, calculations are carried out for different conditions and constraints, depending on the computational size as well as the physical properties and structures of the materials.
Since the different conditions and constraints change coefficient matrices, one needs to solve a large number of different linear systems.
Hence, to use the ILU-based methods listed above, one has no choice but to set up the same parameter value predetermined empirically for each linear system.
However, since every linear system has different optimum parameter values, these methods cannot deliver high performance and may sometimes produce solutions that diverge.
In summary, the previous ILU-based methods have traded in the parameter-free advantage of ILU(0) preconditioning to improve performance, and as a result they have significant problems in practice.

In recent years, new indices to evaluate the effect of ILU preconditioning have been proposed \cite{Duff, Fujiwara, Iwashita2007, Iwashita2005}.
These indices are some functions of ILU-factorized matrices.
However, they cannot be described as explicit functions of the parameters because the parameters directly change not only the values of entries but also the sequence of operations in ILU factorization.
Therefore, they could be available only for brute force optimization as well as the previous ILU-based methods described above.

More sophisticated methods have been proposed to determine the ordering of the ILU factorization using some indices calculated from the matrices being processed \cite{DAZEVEDO}.
The indices are called the minimum update matrix and minimum discarded fill.
Because these self-ordering processes take a tremendously long time, the total calculation time is much increased in almost all cases.
Therefore, these methods are not practical except for the special usage as described above.

With this background, in this paper we attempt to improve ILU preconditioning performance without losing practicality for physical simulations \cite{A2ICCG}.
First, we introduce to ILU preconditioning new acceleration parameters which make themselves easy to be optimized automatically.
We then describe a mechanism which automatically optimizes these acceleration parameters and propose the process as auto-accelerated ILU preconditioning, or A$^2$ILU.
A$^2$ILU can also be applied to the previous ILU-based methods, including shifted ILU, MILU, fill-in level ILU, and Crout ILU.

The rest of the paper is organized as follows.
Section 2 gives an overview of ILU preconditioning.
Section 3 explains the basics of A$^2$ILU preconditioning; acceleration parameters for ILU preconditioning are introduced first, followed by an explanation of the mechanism that automatically optimizes the acceleration parameters.
In section 4, the performance of A$^2$ILU preconditioning is evaluated by numerical experiments.
In subsection 4.1, we evaluate the performance of A$^2$ILU(0) preconditioning with respect to linear systems obtained by discretization of the PDE on rectangular grids and validate its generality, practicality, and scalability.
We also consider applying auto-acceleration to the major ILU-based methods and validate the results.
In subsection 4.2, we evaluate the performance of shifted A$^2$ILU(0) (shifted ILU(0) with auto-acceleration) preconditioning for general sparse matrices by using more than 200 sample matrices obtained from the University of Florida sparse matrix collection \cite{UFSMC}.

\section{ILU preconditioning for sparse matrices}
In this section we present a brief introduction of ILU preconditioning before proposing our method.

Consider a system of linear equations whose coefficients are given by a sparse matrix $A$, denoted as follows:
\begin{equation}
Ax=b.
\label{eqn:linearSystem}
\end{equation}
When such a system is solved by an iterative method, the convergence strongly depends on the property of the coefficient matrix.
It can be expected that the number of iterations decreases as a coefficient matrix tends to be close to an identity matrix.
Preconditioning reconstructs the system as (\ref{eqn:re-constructedLinearSystem}),
\begin{equation}
(K_1^{-1}AK_2^{-1})(K_2x)=K_1^{-1}b,
\label{eqn:re-constructedLinearSystem}
\end{equation}
\begin{equation}
M=K_1K_2.
\label{eqn:preconditioningMatrix}
\end{equation}
Here, the matrix $M$ in (\ref{eqn:preconditioningMatrix}) is called the preconditioning matrix.
The closer $M$ is to $A$, the closer $K_1^{-1}AK_2^{-1}$ (the coefficient matrix after the reconstruction) is to the identity matrix, thus improving the convergence significantly.
ILU preconditioning performs LU factorization on coefficient matrix $A$ in an incomplete way and uses the result as a preconditioning matrix.
LU factorization completely decomposes coefficient matrix $A$ into the product of a strictly lower-triangular matrix $\bar{L}$, a diagonal matrix $\bar{D}$, and a strictly upper-triangular matrix $\bar{U}$,
\begin{equation}
A=(\bar{L}+\bar{D})\bar{D}^{-1}(\bar{D}+\bar{U}).
\label{eqn:LUfactorization}
\end{equation}
On the other hand, ILU factorization keeps a certain degree of sparsity in these matrices by discarding part of the fill-in in the course of the factorization process.
Some sophisticated ILUs such as Crout ILU discard not only fill-in but also the updated original entries.
In particular, ILU(0) requires that these matrices inherit the nonzero entry structure of coefficient matrix $A$ by discarding all fill-ins.
Let $L$ and $U$ denote the strictly lower-triangular matrix and the strictly upper-triangular matrix thus obtained, respectively.
Using these matrices, the ILU preconditioning matrix $M$ can be expressed as follows:
\begin{equation}
M=(L+D)D^{-1}(D+U).
\label{eqn:ILU-preconditioningMatrix}
\end{equation}
Neither $L$ and $\bar{L}$ nor $U$ and $\bar{U}$ generally coincide, so there is a difference between matrices $A$ and $M$.
The difference between $A$ and $M$ is called remainder matrix $R$ as
\begin{equation}
R=A-M,
\label{eqn:remainderMatrix}
\end{equation}
which is used later to evaluate the approximation accuracy of a preconditioning matrix.

\section{Auto-accelerated ILU preconditioner}
In this section we propose the auto-acceleration of ILU preconditioning.
Our new acceleration parameters are used to change only nonzero entry values after ILU factorization.
These parameters can be optimized automatically because the ILU factorized matrix can be described as an explicit function of them.

\subsection{Introduction of acceleration parameters}
In order to improve the computational performance, ILU preconditioning is requested to reduce both the number of iterations and the operation count per iteration step.
The former can be obtained by the reconstruction of a coefficient matrix shown in (\ref{eqn:re-constructedLinearSystem}).
The latter can be given by keeping the sparsity of a coefficient matrix.
However, the sparsity degrades the approximation accuracy of the preconditioned matrix resulting in an increase in the number of iterations.

To achieve a balance between the two factors, ILU preconditioning performs LU factorization on coefficient matrix $A$ under some sort of constraints related to the structure of the nonzero entries.
However, the effects of the constraints on the accuracy of ILU factorization are quantitatively unknown.
There is no guarantee that the preconditioning matrix obtained by ILU factorization has an optimum level of approximation for coefficient matrix $A$ while preserving its nonzero entry structure.
We focus on this point and attempt to improve the approximation accuracy of coefficient matrix $A$ by modifying only the values (without changing the nonzero entry structure) of the preconditioning matrix obtained by ILU factorization.
If such an attempt is successful, the number of iterations can be reduced without increasing the computational time for each iteration process, ensuring faster computation without any trade-off.
In our proposed method, we introduce distinct scalar parameters for the matrices obtained by ILU(0) or the other ILU-based methods: $\phi$ for strictly triangular matrices $L$ and $U$, and $\gamma$ for diagonal matrix $D$,
\begin{equation}
M=(\phi L+\gamma D)(\gamma D)^{-1}(\gamma D+\phi U).
\label{eqn:A2ILU-preconditioningMatrix}
\end{equation}
In the rest of this paper, we refer to $\phi$ and $\gamma$ as acceleration parameters.

\subsection{Automatic optimization of the acceleration parameters}
We discuss the mechanism which automatically optimizes the acceleration parameters.
From the discussion in the previous section, minimizing the approximation error of the preconditioning matrix relative to the coefficient matrix is thought to be equal to improving the convergence of the iterative method maximally.
Therefore, if we express the approximation error with some kind of objective function of $R$, we can optimize the acceleration parameters automatically with gradient-based methods because the objective function of $R$ can be written explicitly as a function of these acceleration parameters.

We adopt the squared Euclidean norm of $Re$, where $e = (1,\ldots, 1)^T$ as the objective function of $R$,
\begin{equation}
f(R)=\|Re\|_2^2=\left\|\sum_{j=0}^{n-1}r_{ij}\right\|_2^2=\sum_{i=0}^{n-1}\left(\sum_{j=0}^{n-1}r_{ij}\right)^2.
\label{eqn:objectiveFunction}
\end{equation}
Here, $n$ is the size of the matrix, $i$ is the index for the rows, and $j$ is the index for the columns.
The objective function is based on an idea of MILU in which the condition number, and consequently the number of iterations, is reduced by minimizing $(A-M)e = Re$ for solving an elliptic PDE.

The objective function $f(R)$ is written as a nonlinear explicit function of $\phi$ and $\gamma$,
\begin{eqnarray}
f(R)=\sum_{i=0}^{n-1}\left(\sum_{j=0}^{n-1}\left(a_{ij}-\phi l_{ij}-\gamma d_{ij}-\phi u_{ij}-\phi^2\gamma^{-1}\sum_{k=0}^{\min[i,j]-1}l_{ik}{d_{kk}}^{-1}u_{kj}\right)\right)^2. \nonumber \\
\label{eqn:objectiveFunctionExpanded}
\end{eqnarray}

The optimum values of the acceleration parameters are obtained by minimizing (\ref{eqn:objectiveFunctionExpanded}).
When we use Newton--Raphson method to optimize, an update equation for the parameters is given as
\begin{equation}
\left[\begin{array}{c}\phi^{(t+1)} \\ \gamma^{(t+1)} \end{array}\right]=
\left[\begin{array}{c}\phi^{(t)} \\ \gamma^{(t)} \end{array}\right]-H^{-1}g.
\label{eqn:newtonRaphsonMethod}
\end{equation}
Here, $H$ and $g$ are, respectively, a Hessian matrix and a gradient vector defined by
\begin{eqnarray}
H&=&\left[\begin{array}{cc}(\partial^2/\partial\phi^2)f(R)&(\partial^2/\partial\phi\partial\gamma)f(R) \\(\partial^2/\partial\gamma\partial\phi)f(R)&(\partial^2/\partial\gamma^2)f(R)\end{array}\right],
\label{eqn:hesseMatrix} \\
g&=&\left[\begin{array}{cc}((\partial/\partial\phi)f(R)&(\partial/\partial\gamma)f(R)\end{array}\right]^T.
\label{eqn:gradientVector}
\end{eqnarray}
A pair of the acceleration parameters $\phi$ and $\gamma$ in ({\ref{eqn:A2ILU-preconditioningMatrix}}) and the optimization of them by minimizing the $f(R)$ in ({\ref{eqn:objectiveFunction}}) are referred to as auto-acceleration.
Since this auto-acceleration is applied to a matrix described in (\ref{eqn:ILU-preconditioningMatrix}), it is expected that the auto-acceleration is incorporated with any other ILU-based methods.
Applications of the auto-acceleration to major ILU-based methods are validated by numerical experiments as well as ILU(0).
In the rest of this paper, we abbreviate ILU-based preconditioning method with auto-acceleration to A$^2$ILU.
For example, shifted ILU(0) with auto-acceleration is referred to as shifted A$^2$ILU(0).

\subsection{Computational cost of auto-acceleration}
Here we evaluate the computational cost of auto-acceleration by comparing it with MILU factorization.
It can be helpful for the evaluation to divide the computational cost into three parts, i.e., costs of loading from memory, operating on CPU, and storing into memory.
Hereafter, we refer to these three costs as loading cost, operating cost, and storing cost, respectively.

In the auto-acceleration, the matrix-matrix product $l_{ik}d_{kk}^{-1}u_{kj}$ in (\ref{eqn:objectiveFunctionExpanded}) dominates its computational cost.
On the other hand, MILU performs the following operations:
\begin{eqnarray}
l_{ij}+d_{ij}+u_{ij} & = &
a_{ij}-\sum_{k=0}^{\min[i,j]-1}l_{ik}{d_{kk}}^{-1}u_{kj} \mbox{ if} \ (i, j) \in P {\rm ,} 
\label{eqn:ILUfactorization} \\
d_{ii} & = & d_{ii} - \omega \sum_{k=0}^{\min[i,j]-1}l_{ik}d_{kk}^{-1}u_{kj} \mbox{ otherwise},
\label{eqn:MILUfactorization}
\end{eqnarray}
where $P$ is a set of indices with nonzero entries of the resultant preconditioned matrix.
MILU updates $L$, $D$, and $U$ by using (\ref{eqn:ILUfactorization}) with an additional modification in (\ref{eqn:MILUfactorization}).
We can see that both equations are also dominated by the matrix-matrix product.
Since the matrix-matrix product must be done for all nonzero elements in both methods, their loading and operating costs are, respectively, identical.

As for the storing cost, we can see a difference between these methods.
The auto-acceleration stores only two scalar acceleration parameters, $\phi$ and $\gamma$, and a constant number of additional variables temporarily used in (\ref{eqn:newtonRaphsonMethod}) to (\ref{eqn:gradientVector}), while MILU stores all nonzero elements of three matrices, $L$, $D$, and $U$.
Therefore, the storing cost of MILU is much larger than that of the auto-acceleration for large practical problems.

The computational cost of a simple ILU is less than MILU by the matrix-matrix product of the additional modification in (\ref{eqn:MILUfactorization}).
The comparison between the auto-acceleration and ILU is much more complicated than the case of the auto-acceleration and MILU.
Computational cost of auto-acceleration in sample problems will be compared with that of whole A$^2$ILU preconditioned iterative method in later sections.

\subsection{Constraint on the acceleration parameters}
By preliminary numerical experiments, we found that the parameters satisfying $\gamma / \phi > 1.0$ degraded the computational performance in almost all cases.
Although the theoretical aspect of this problem is under investigation, we employ the condition $\gamma / \phi \leq 1.0$ as an empirical constraint on the parameter values.
In numerical experiments described in the subsequent section, the acceleration parameters of all auto-accelerated methods are subject to this constraint.

\section{Numerical experiments}
The proposed method is validated by numerical experiments.
Numerical experiments described below are classified into two groups according to the type of coefficient matrix used.
The first one uses multidiagonal sparse matrices given by discretizing partial differential equations (PDEs) on rectangular grids, while the second one uses general sparse matrices obtained from the University of Florida sparse matrix collection.
Practicality and performance of A$^2$ILU and auto-accelerated ILU-based methods are evaluated in the first experiment.
The second experiment mainly shows the effect of the auto-acceleration on shifted ILU(0).

\subsection{Performance evaluation for linear systems arising from rectangular grids}
In this subsection we evaluate the performance of A$^2$ILU preconditioning for systems of linear equations obtained through discretization of PDEs on rectangular grids.
We also validate its effectiveness from various aspects.
Here, the coefficient matrices of the system of equations will be multidiagonal (e.g., tridiagonal, pentadiagonal) matrices where multiple arrays of nonzero entries are lined up on and along the main diagonal.

\subsubsection{Generality for various types of physical simulations}
Here we validate the generality of A$^2$ILU(0) preconditioning using linear systems derived from five different types of physical simulation.
Table \ref{table:physicalSimulations} shows the details of each physical simulation.
Type of PDE and Stationary refer to the specifications of PDEs obtained by modeling physical phenomena.
These PDEs were discretized on rectangular grids using the method described in the Discretization column to produce linear systems.
Table \ref{table:matricesAndSolver1} shows the details of the coefficient matrix of the linear system and the iterative method used.
The size of the coefficient matrix is denoted by n, and nnz/row refers to the number of nonzero entries for each row.
The convergence criterion of the iterative methods is given by $\|r\|^2/\|b\|^2 \le \epsilon$, where $r$ is the recursive residual vector calculated by the recurrence formula, $r_{k+1}=r_k-\alpha_kAp_k$ with the CG method.
The threshold value $\epsilon$ is determined to avoid pseudo-convergence so that $\|s\|^2/\|r\|^2 \le 2$ is satisfied until the iterative method is terminated, where $s$ is the true residual vector calculated by $s_k=b-Ax_k$.
In each problem we carry out diagonal scaling first to normalize the coefficient matrix so that every diagonal entry is 1.
Table \ref{table:calculationEnvironment1} summarizes the details of the computation environment used.

\begin{table}[!h]
\caption{Physical simulations examined.} 
\begin{center} \footnotesize
\begin{tabular}{|c|c|c|c|c|} \hline
No. & Physics & Type of PDE & Stationary & Discretization \\ \hline
1 & Incompressible fluid & Poisson & Yes & FVM \\ \hline
2 & Heat conduction & Poisson & No & FVM \\ \hline
3 & Light diffusion & Helmholtz & Yes & FEM \\ \hline
4 & Heat radiation & Helmholtz & No & FVM \\ \hline
5 & Charge transfer & Advection-diffusion & No & FVM \\ \hline
\end{tabular}
\end{center}
\label{table:physicalSimulations}
\end{table}

\begin{table}[!h]
\caption{Properties of coefficient matrices and iterative method used.}
\begin{center} \footnotesize
\begin{tabular}{|c|c|c|c|c|c|} \hline
No. & $n$ & nnz / row & Symmetry & Solver & $\epsilon$ \\ \hline
1 & 537600 & 7 & Yes & CG \cite{CG} & 1.0e$-$11 \\ \hline
2 & 720000 & 7 & Yes & CG & 1.0e$-$12 \\ \hline
3 & 410913 & 27 & Yes & CG & 1.0e$-$12 \\ \hline
4 & 647168 & 7 & Yes & CG & 1.0e$-$11 \\ \hline
5 & 884736 & 7 & No & BiCGSTAB \cite{BiCGSTAB} & 1.0e$-$9 \\ \hline
\end{tabular}
\end{center}
\label{table:matricesAndSolver1}
\end{table}

\begin{table}[!h]
\caption{Calculation environment.}
\begin{center} \footnotesize
\begin{tabular}{|c|c|c|} \hline
System & NEC LX118Tc-4G \\ \hline
\multirow{2}{*}{CPU} & Intel Xeon X5670 \\
& 2.93 GHz / L3-Cache 12 MB \\ \hline
Memory & DDR3-1333 48 GB \\ \hline
Compiler & Intel C++ Ver.11.1.075 \\ \hline
\end{tabular}
\end{center}
\label{table:calculationEnvironment1}
\end{table}

Each problem was solved using both ILU(0) preconditioning and A$^2$ILU(0) preconditioning.
Table \ref{table:result1} shows the results of solving the linear system involved in each physical simulation.
Itr. denotes the number of iterations required for convergence.
Time denotes the number of seconds required for convergence; in A$^2$ILU(0) the time for the auto-acceleration is shown in brackets.
$\phi$ and $\gamma$ in A$^2$ILU(0) are the values optimized by the auto-acceleration process.
Speed-up ratio is the quotient obtained when the number of iterations or the computational time for ILU(0) is divided by the corresponding value for A$^2$ILU(0); they show how much the performance is improved by auto-acceleration.

\begin{table}[!h]
\caption{Performance of ILU{\rm (0)} and A$^2$ILU{\rm (0)} for sample problems.}
\begin{center} \footnotesize
\begin{tabular}{|c|c|c|c|c|c|c|c|c|} \hline
\multirow{3}{*}{No.} & \multicolumn{2}{|c|}{\multirow{2}{*}{ILU(0)}} & \multicolumn{4}{|c|}{\multirow{2}{*}{A$^2$ILU(0)}} & \multicolumn{2}{|c|}{Speed-up} \\
& \multicolumn{2}{|c|}{} & \multicolumn{4}{|c|}{} & \multicolumn{2}{|c|}{ratio} \\ \cline{2-9}
& Itr. & Time & Itr. & Time & $\phi$ & $\gamma$ & Itr. & Time \\ \hline
1 & 168 & 4.34 & 96 & 2.59 (2.53e$-$2) & 1.02 & 0.66 & 1.75 & 1.68 \\ \hline
2 & 125 & 6.50 & 72 & 3.95 (5.10e$-$2) & 1.52 & 0.98 & 1.74 & 1.65 \\ \hline
3 & 81 & 6.94 & 46 & 4.76 (8.88e$-$2) & 1.13 & 0.80 & 1.76 & 1.46 \\ \hline
4 & 191 & 8.41 & 97 & 4.49 (4.46e$-$2) & 1.53 & 0.98 & 1.97 & 1.87 \\ \hline
5 & 93 & 10.51 & 55 & 6.50 (7.18e$-$2) & 0.80 & 0.51 & 1.69 & 1.62 \\ \hline
\end{tabular}
\end{center}
\label{table:result1}
\end{table}

According to Table \ref{table:result1}, performance was improved for every problem when auto-acceleration was used.
The average speed-up ratio for iterations was 1.78 with the range of [1.69,1.97].
The average speed-up ratio for time was 1.65 with the range of [1.46,1.87].
Meanwhile, acceleration parameters $\phi$ and $\gamma$ were optimized to different values for each problem.
The optimum values of these acceleration parameters differ depending on the physical phenomena involved and the type of PDEs used.
The time of the auto-acceleration is just 1\% to 2\% of the total computing time of A$^2$ILU(0).

\subsubsection{Comparison of A$^2$ILU(0) and ILU-based methods}
In this subsection we compare A$^2$ILU(0) preconditioning with the various ILU-based methods that have been previously proposed and validate the superiority of A$^2$ILU(0) preconditioning.
We chose the following four as major valid methods:

\begin{enumerate}

\item Shifted ILU.
ILU factorization in this method is performed for the matrix $\tilde{A}$ obtained by shifting the diagonal entries of coefficient matrix $A$.
ILU corresponds to the shift parameter $\alpha=0$.
We use the shifted ILU(0)
\begin{equation}
\tilde{A}=A+\alpha\mbox{diag}(A).
\label{eqn:shiftedILU}
\end{equation}

\item MILU.
In this ILU factorization, the fill-in being discarded is multiplied by a relaxing parameter $\omega$, and the product is subtracted from the diagonal entry of that row.
ILU corresponds to the case where $\omega=0$.
We use the MILU(0)
\begin{equation}
d_{ii}=d_{ii}-\omega \sum_{k=0}^{\min[i,j]-1}l_{ik}{d_{kk}}^{-1}u_{kj} \qquad {\rm if} \ i \ne j \ {\rm and} \ a_{ij}=0.
\label{eqn:modifiedILU}
\end{equation}

\item Fill-in level ILU.
In this ILU factorization, fill-ins are allowed as long as the fill-in level is below or equal to $p$.
The fill-in level at entry $(i, j)$ is updated during ILU factorization as the nonzero entries are updated.
The initial value of the fill-in level is given by (\ref{eqn:fill-inLevelILU1}) and is updated by (\ref{eqn:fill-inLevelILU2}).
When $p=0$, it amounts to the ILU(0)
\begin{eqnarray}
lev_{ij}&=&\left \{ \begin{array}{ll}
0 & {\rm if} \ a_{ij}=0 \ {\rm or} \ i=j, \\
\infty & {\rm otherwise} ,
\end{array} \right.
\label{eqn:fill-inLevelILU1} \\
lev_{ij}&=&{\rm min}\{lev_{ij}, lev_{ik}+lev_{kj}+1\}.
\label{eqn:fill-inLevelILU2}
\end{eqnarray}

\item Crout ILU.
In this ILU factorization, $L$ is accessed by columns and $U$ by rows in the same way as in Crout LU factorization.
Any fill-in whose impact on $L^{-1}$ or $U^{-1}$ is less than a tolerance is discarded and the number of fill-ins in each column of $L$ or each row of $U$ is limited.
The limit is obtained by $nnz / 2n * m$, where $m$ is a ratio of maximum fill-in.
When the drop tolerance $tol=0$ and $m=\infty$, it amounts to the complete LU factorization.

\end{enumerate}

We have five parameters to be specified before running the methods: shift parameter $\alpha$, relaxing parameter $\omega$, fill-in level $p$, drop tolerance $tol$, and ratio of maximum fill-in $m$.
Since their optimum values are unknown, each method is performed for a set of candidate values of parameters as follows:
\begin{eqnarray*}
\alpha & = & -0.4 + 0.1 j, \quad j \in \{ 0, 1,\ldots, 10 \}, \\
\omega & = & -0.5 + 0.1 j, \quad j \in \{ 0, 1,\ldots, 16 \}, \\
p & = & \{ 1, 2, 3 \}, \\
m & = & \{ 1, 2, 5, 10 \}, \\
tol & = & \{ 0.001, 0.002, 0.004, 0.01, 0.02, 0.04, 0.1, 0.2 \}.
\end{eqnarray*}

\begin{figure}[!h]
\begin{center}
\includegraphics[scale=\figscale,clip]{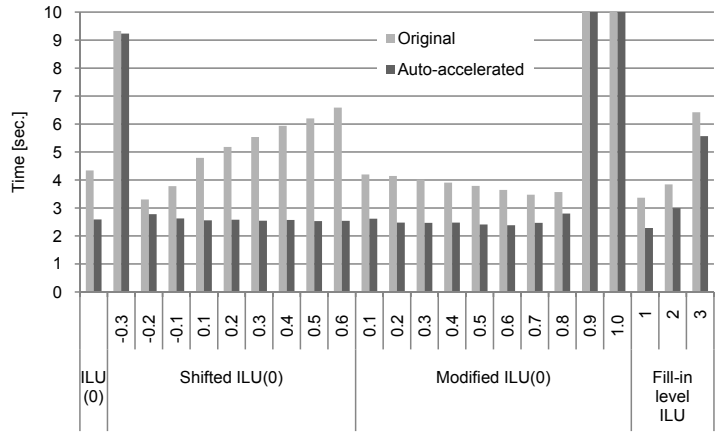}
\caption{Result of sample problem {\rm 1} by all methods except Crout ILU.}
\label{fig:resultOfILUs1}
\end{center}
\end{figure}

\begin{figure}[!h]
\begin{center}
\includegraphics[scale=\figscale,clip]{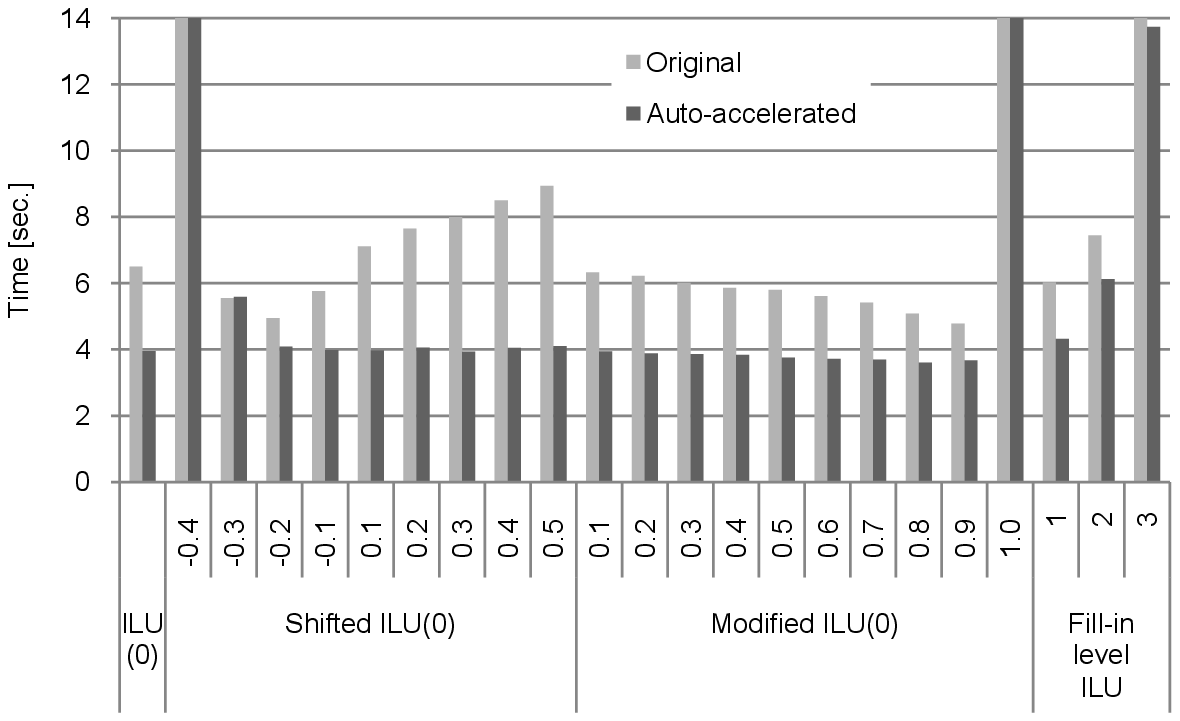}
\caption{Result of sample problem {\rm 2} by all methods except Crout ILU.}
\label{fig:resultOfILUs2}
\end{center}
\end{figure}

\begin{figure}[!h]
\begin{center}
\includegraphics[scale=\figscale,clip]{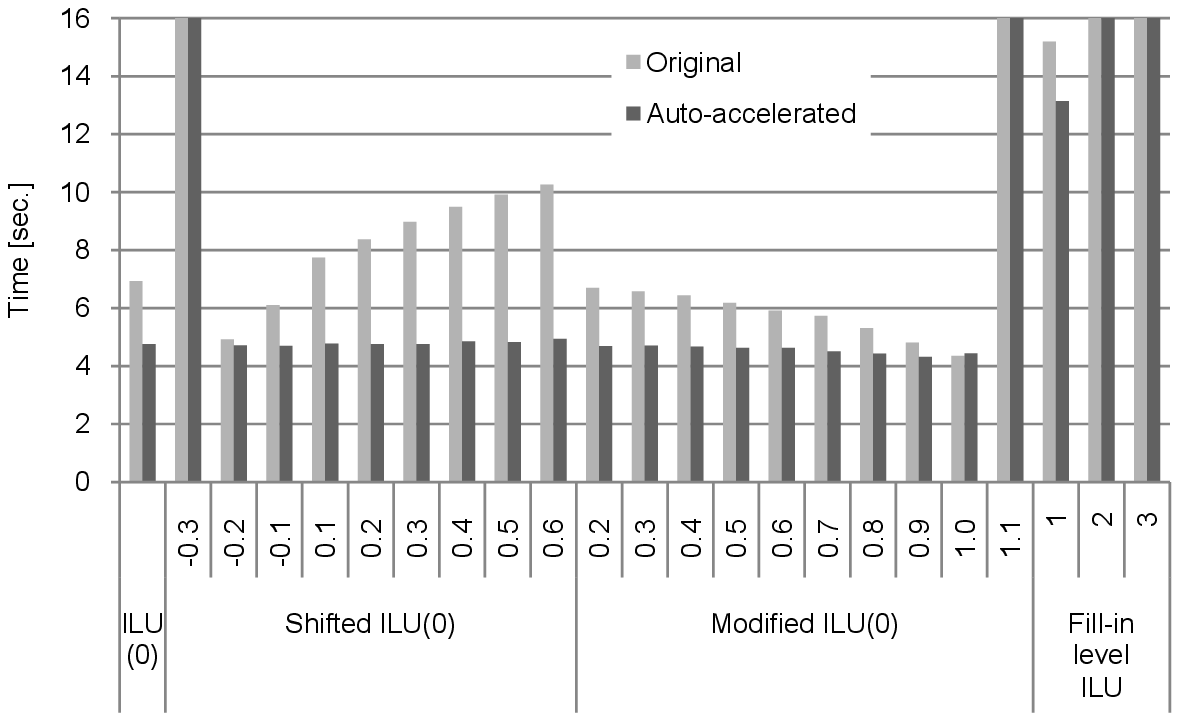}
\caption{Result of sample problem {\rm 3} by all methods except Crout ILU.}
\label{fig:resultOfILUs3}
\end{center}
\end{figure}

\begin{figure}[!h]
\begin{center}
\includegraphics[scale=\figscale,clip]{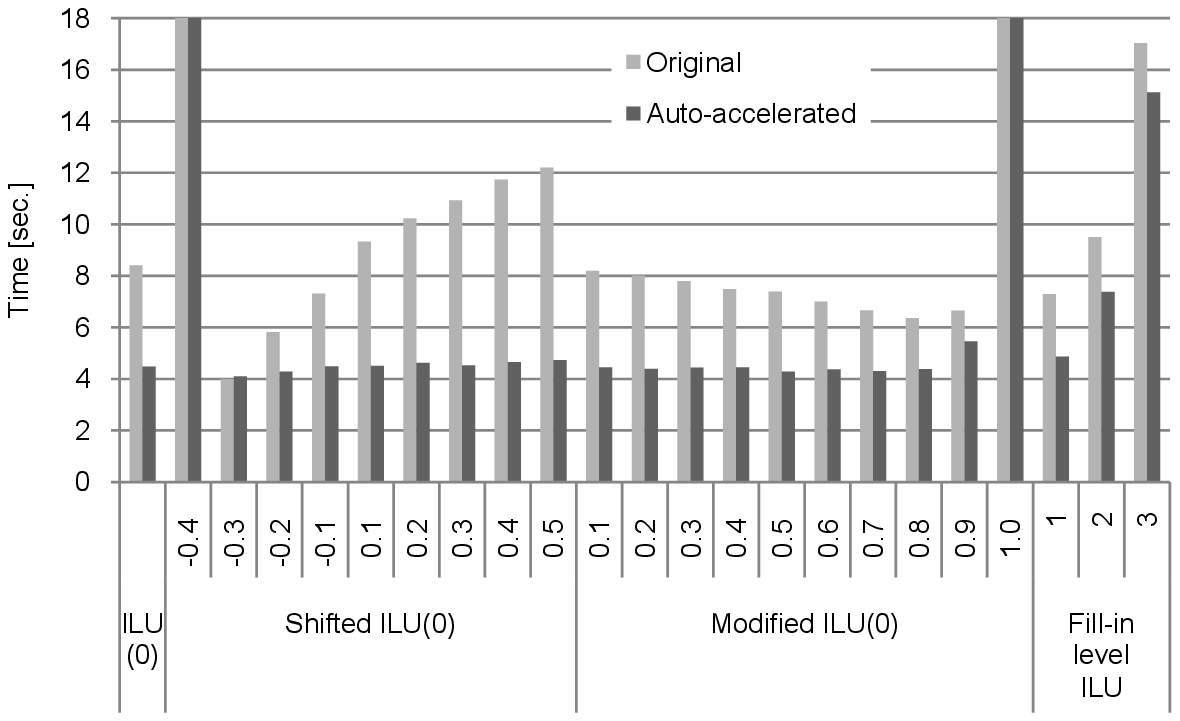}
\caption{Result of sample problem {\rm 4} by all methods except Crout ILU.}
\label{fig:resultOfILUs4}
\end{center}
\end{figure}

\begin{figure}[!h]
\begin{center}
\includegraphics[scale=\figscale,clip]{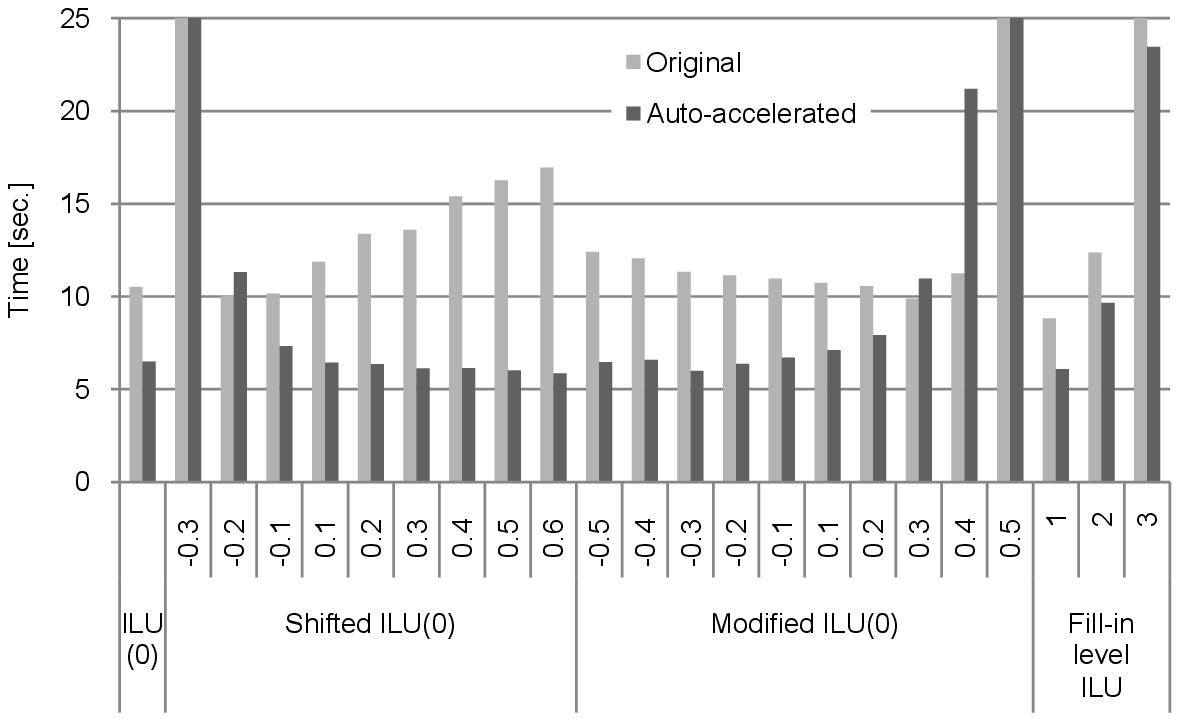}
\caption{Result of sample problem {\rm 5} by all methods except Crout ILU.}
\label{fig:resultOfILUs5}
\end{center}
\end{figure}

We used the same sample problems and computation environment as in subsection 4.1.1.
Tables \ref{table:physicalSimulations}, \ref{table:matricesAndSolver1}, and \ref{table:calculationEnvironment1} show the details.
Figures \ref{fig:resultOfILUs1} through \ref{fig:resultOfILUs5} show the results of sample problems 1 through 5, respectively, by all methods except Crout ILU.
From the left side along the horizontal axis, the graph shows ILU(0), shifted ILU(0), MILU(0), and fill-in level ILU, in that order.
The values on the horizontal axis indicate parameters $\alpha$, $\omega$, and $p$ for the respective methods.
The vertical axis represents the computational time required for convergence (in seconds).
The light and dark gray bars show the results of previous methods and auto-accelerated previous methods, respectively.
The leftmost light and dark gray bars are, respectively, ILU(0) and A$^2$ILU(0).

\begin{figure}[!h]
\begin{center}
\includegraphics[scale=\figscale,clip]{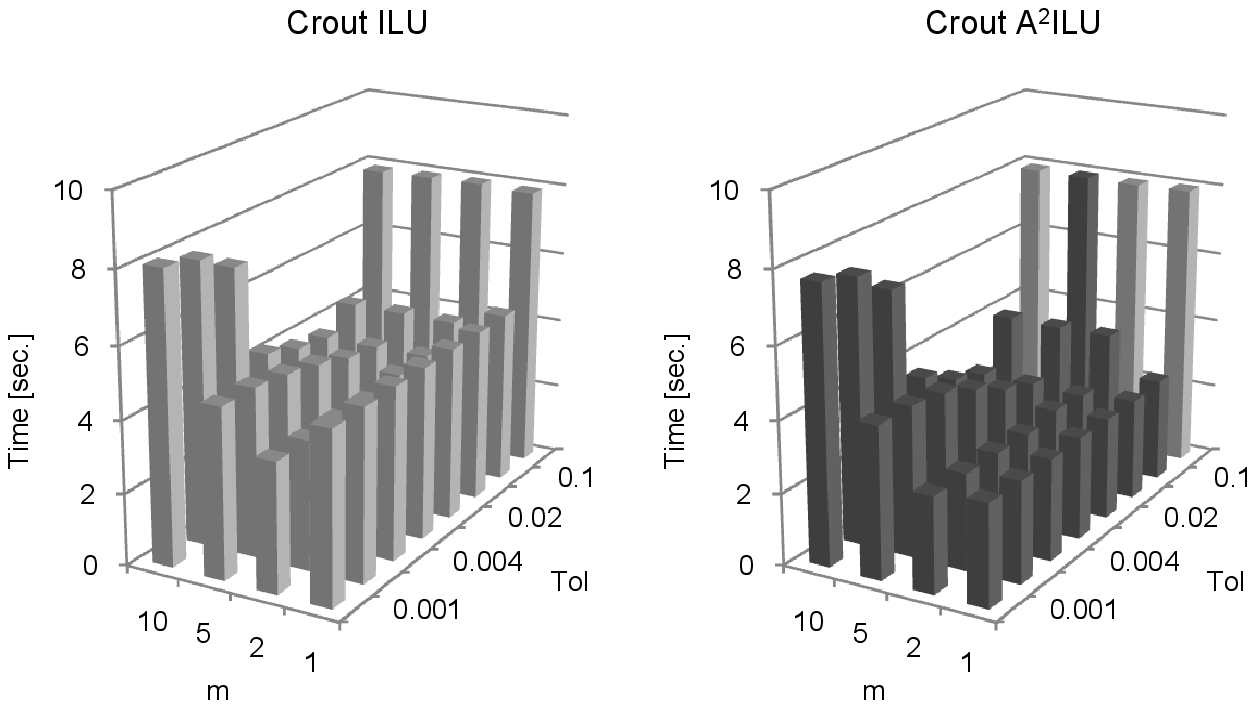}
\caption{Result of sample problem {\rm 1} by Crout ILU. The dark gray bars in the results of Crout A$^2$ILU denote the cases when the auto-acceleration reduces the computing time of Crout ILU.}
\label{fig:resultOfCroutILU1}
\end{center}
\end{figure}

\begin{figure}[!h]
\begin{center}
\includegraphics[scale=\figscale,clip]{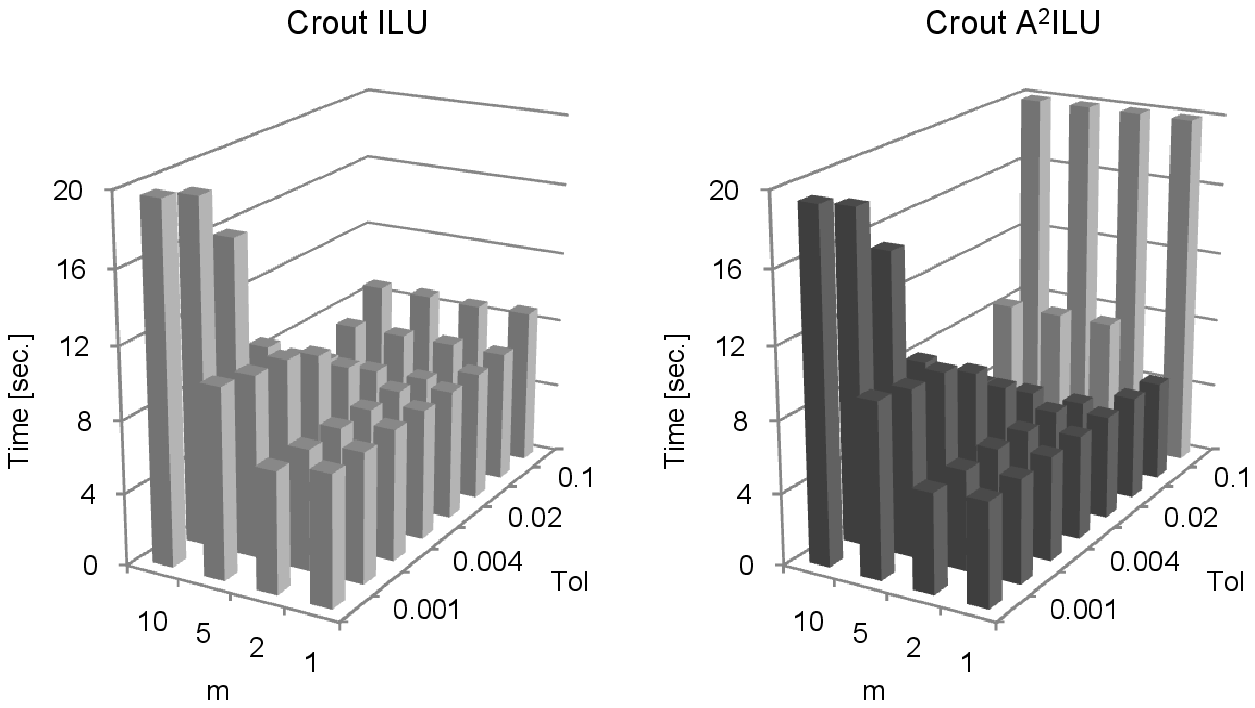}
\caption{Result of sample problem {\rm 2} by Crout ILU. The dark gray bars in the results of Crout A$^2$ILU denote the cases when the auto-acceleration reduces the computing time of Crout ILU.}
\label{fig:resultOfCroutILU2}
\end{center}
\end{figure}

\begin{figure}[!h]
\begin{center}
\includegraphics[scale=\figscale,clip]{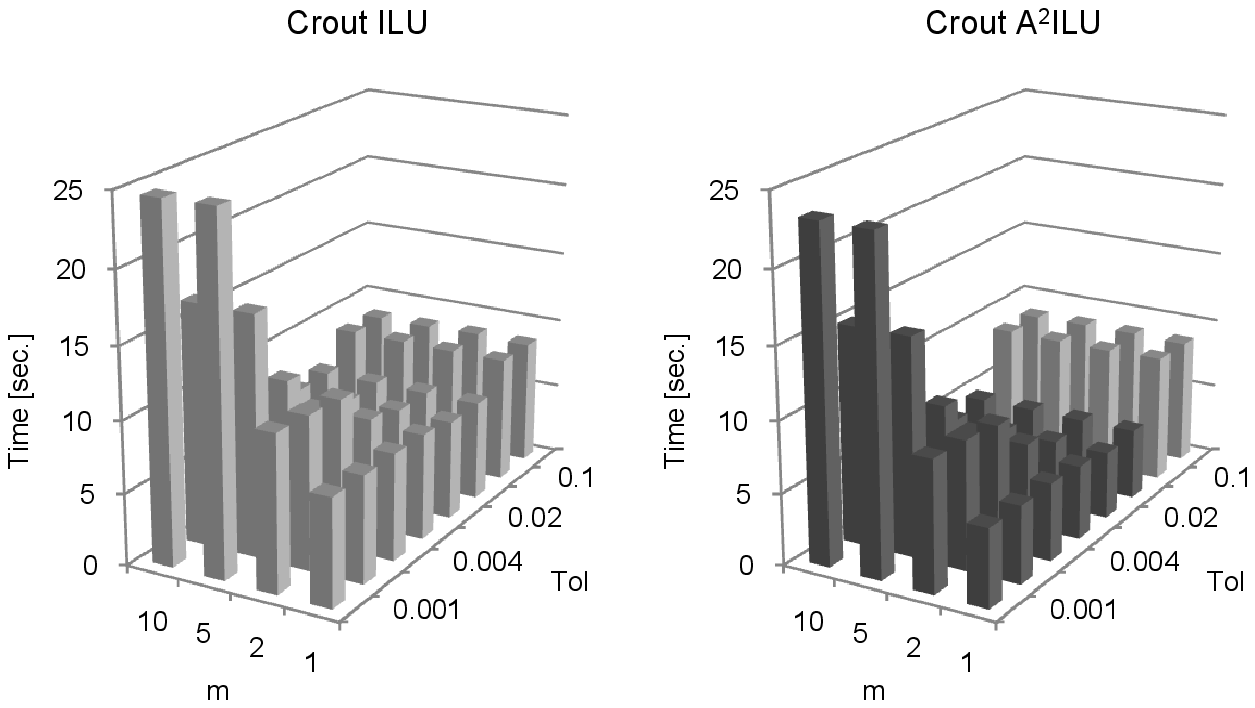}
\caption{Result of sample problem {\rm 3} by Crout ILU. The dark gray bars in the results of Crout A$^2$ILU denote the cases when the auto-acceleration reduces the computing time of Crout ILU.}
\label{fig:resultOfCroutILU3}
\end{center}
\end{figure}

\begin{figure}[!h]
\begin{center}
\includegraphics[scale=\figscale,clip]{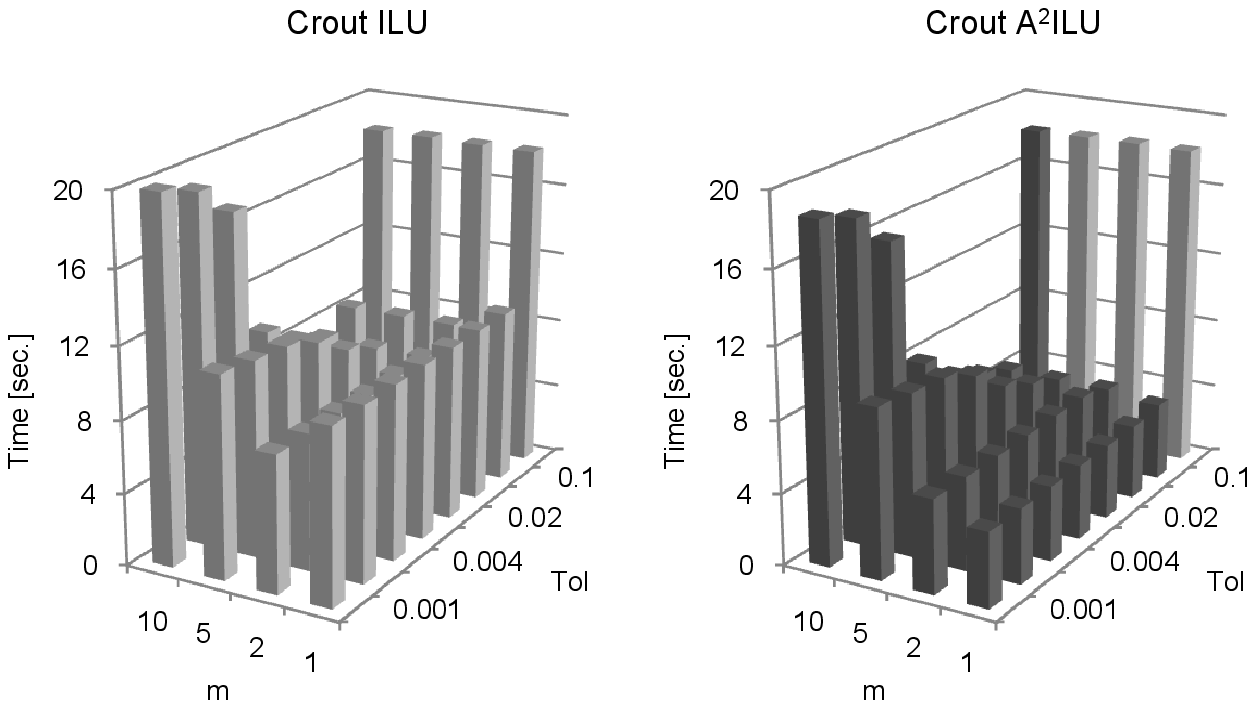}
\caption{Result of sample problem {\rm 4} by Crout ILU. The dark gray bars in the results of Crout A$^2$ILU denote the cases when the auto-acceleration reduces the computing time of Crout ILU.}
\label{fig:resultOfCroutILU4}
\end{center}
\end{figure}

\begin{figure}[!h]
\begin{center}
\includegraphics[scale=\figscale,clip]{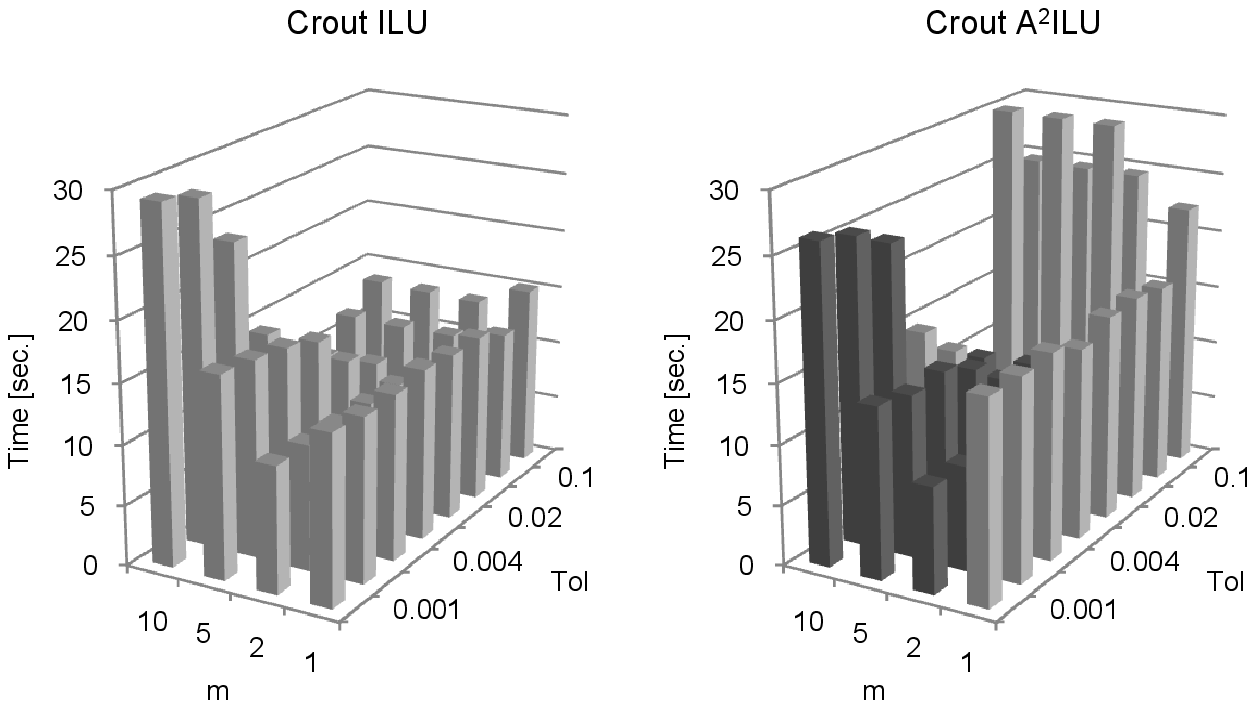}
\caption{Result of sample problem {\rm 5} by Crout ILU. The dark gray bars in the results of Crout A$^2$ILU denote the cases when the auto-acceleration reduces the computing time of Crout ILU.}
\label{fig:resultOfCroutILU5}
\end{center}
\end{figure}

Figures \ref{fig:resultOfCroutILU1} through \ref{fig:resultOfCroutILU5} show the results of sample problems 1 through 5, respectively, by Crout ILU.
The horizontal axis represents drop tolerance $tol$ and ratio of maximum fill-in $m$.
The vertical axis represents the computational time required for convergence in seconds.
The results by Crout ILU and by Crout A$^2$ILU are shown in the left and right graphs, respectively.
The dark gray bars in the results of Crout A$^2$ILU denote the cases when the auto-acceleration reduces the computing time of Crout ILU.

First, we compare these results of A$^2$ILU(0) preconditioning and the previous ILU-based methods from the standpoint of practicality.
In each of the four ILU-based methods, the value of the parameter drastically influences the performance.
See the light gray bars for shifted ILU(0), MILU(0), and fill-in level ILU in Figures \ref{fig:resultOfILUs1} through \ref{fig:resultOfILUs5}.
See the left-hand graphs in Figures \ref{fig:resultOfCroutILU1} through \ref{fig:resultOfCroutILU5} for Crout ILU.
If the parameter is optimized, the computations by all these ILU-based methods take less time to converge than ILU(0).
However, the farther the parameter is from its optimum value, the larger the computational time.
In several cases it takes more time to converge than the ILU(0), and in some cases the residual norm diverges, resulting in no convergent solutions.
To avoid such unacceptable situations, the parameter must be optimized in advance.
However, the optimum value of each parameter differs from problem to problem.
For instance, the optimum value of the relaxing parameter $\omega$ in MILU varies from 0.3 to 1.0 depending on the problem.
The parameters must therefore be optimized individually for each problem.
However, it is essentially impossible to optimize the parameters in such an efficient way that it will not adversely affect the performance improvement made over the ILU(0).
For instance, in Figure \ref{fig:resultOfILUs1}, we obtain the optimum values of $\alpha$ for shifted ILU(0) by brute force searching which takes the sum of the computational times over the set of candidate parameter values, more than 10 times that of ILU(0) preconditioning.
A$^2$ILU(0) preconditioning has none of these problems since the acceleration parameters are automatically optimized in a short time.

Next, we compare A$^2$ILU(0) preconditioning and the previous methods from the standpoint of performance.
According to Figures \ref{fig:resultOfILUs1} to \ref{fig:resultOfCroutILU5}, A$^2$ILU(0) preconditioning takes less time for convergence than any of the previous methods in every trial result (Figures \ref{fig:resultOfILUs1} to \ref{fig:resultOfILUs5}: all light gray bars; Figures \ref{fig:resultOfCroutILU1} to \ref{fig:resultOfCroutILU5}: all left-hand graphs) except for a few cases (MILU(0) with $\omega = 1.0$ in Figure \ref{fig:resultOfILUs3}, shifted ILU(0) with $\alpha = -0.3$ in Figure \ref{fig:resultOfILUs4}).
This means that A$^2$ILU(0) preconditioning is better than almost all cases even if the optimum value of the parameter for each previous method can be predetermined by any means.

The above results show that A$^2$ILU(0) preconditioning is superior to the previous ILU-based methods from the standpoints of both performance and practicality for multidiagonal matrices obtained from PDEs.

\subsubsection{Application of auto-acceleration to previous ILU-based methods}
In this subsection we validate the effects when the auto-acceleration is incorporated into the previous ILU-based methods.
In Figures \ref{fig:resultOfILUs1} through \ref{fig:resultOfILUs5} the dark gray bars show the results when the auto-acceleration was applied to shifted ILU(0), MILU(0), and fill-in level ILU.
The auto-acceleration for these methods decreases computational time in every problem regardless of the values of their original parameters.

When the auto-acceleration was applied to shifted ILU(0) and MILU(0), the parameter responsiveness became rather consistent, except for extremely expensive computational time-consuming cases where some diagonal elements approach zero too closely to be improved by auto-acceleration.
Note that MILU(0) with $\omega = 1$ minimizes the objective function $f(R) = 0$ so that the auto-acceleration is fruitless.
For most parameter values, the convergence took less time than ILU(0).
Hence, even without optimizing their original parameters for each problem, the performance was improved over ILU(0).
This suggests that the auto-acceleration improves the practicality of shifted ILU(0) and MILU(0) as well.

Figures \ref{fig:resultOfCroutILU1} through \ref{fig:resultOfCroutILU5} show the effect of auto-acceleration on Crout ILU.
The left and the right figures depict the results before and after application of the auto-acceleration, respectively.
For all sample problems shown in Figures \ref{fig:resultOfCroutILU1} through \ref{fig:resultOfCroutILU5},
the auto-acceleration decreases the computing time for greater $m$ and less $tol$, while the smallest $m$ or greater $tol$ gives worse results.
Considering the roles of $m$ and $tol$, it can be expected that a parameter set of greater $m$ and less $tol$ makes a preconditioned matrix denser.

\subsubsection{Scalability for the coefficient matrix size}
In this subsection we validate the scalability of A$^2$ILU(0) preconditioning with respect to the number of unknowns in the system of linear equations, i.e., the size of the coefficient matrix.
We used a Dirichlet boundary-value problem of a three-dimensional Poisson equation as a sample problem,
\begin{eqnarray}
& -\nabla \cdot (\kappa \nabla u)=f \\
& {\rm in} \ \Omega=(0,1)^3, \nonumber \\
& u(x,y,z)=0 \ {\rm on} \ \partial \Omega , \nonumber \\
& \kappa (x,y,z)=\left \{ \begin{array}{ll}
10^3 & {\rm if} \ {{1}\over{4}} \le x,y,z \le {{3}\over{4}}, \\
1 & {\rm otherwise}, \\
\end{array} \right. \nonumber \\
& f(x,y,z)=x+y+z. \nonumber
\label{eqn:sampleProblem}
\end{eqnarray}

For discretization of the equations, we used a rectangular grid described as before.
The solution vector $u$ was initialized to be 0.
The convergence criterion of the iterative methods is given by $\|r\|^2/\|b\|^2 \le \epsilon$.
The diagonal scaling was carried out in advance.
The number of lattice points on each axis was set to 10 first and intermittently changed up to 640.
The specifications of the coefficient matrix and the iterative method are shown in Table \ref{table:matricesAndSolver2}.
Details of the computation environment used are shown in Table \ref{table:calculationEnvironment1}.
The results are shown in Tables \ref{table:result2} and \ref{table:result3}.

\begin{table}[!h]
\caption{Properties of the matrices and the method used.}
\begin{center} \footnotesize
\begin{tabular}{|c|c|c|c|c|} \hline
$n$ & nnz / row & Symmetry & Solver & $\epsilon$ \\ \hline
from $10^3$ to $640^3$ & 7 & Yes & CG & 1.0e$-$9 \\ \hline
\end{tabular}
\end{center}
\label{table:matricesAndSolver2}
\end{table}

\begin{table}[!h]
\caption{Scalability of the speed-up ratio by A$^2$ILU{\rm (0)}.}
\begin{center} \footnotesize
\begin{tabular}{|c|c|c|c|c|c|c|c|c|} \hline
\multirow{3}{*}{$n$} & \multicolumn{2}{|c|}{\multirow{2}{*}{ILU(0)}} & \multicolumn{4}{|c|}{\multirow{2}{*}{A$^2$ILU(0)}} & \multicolumn{2}{|c|}{Speed-up} \\
& \multicolumn{2}{|c|}{} & \multicolumn{4}{|c|}{} & \multicolumn{2}{|c|}{ratio} \\ \cline{2-9}
& Itr. & Time & Itr. & Time & $\phi$ & $\gamma$ & Itr. & Time \\ \hline
$10^3$ & 21 & 1.87e$-$3 & 19 & 1.99e$-$3 & 1.38 & 1.03 & 0.95 & 0.94 \\ \hline
$20^3$ & 33 & 1.95e$-$2 & 27 & 1.70e$-$2 & 1.86 & 1.24 & 1.22 & 1.15 \\ \hline
$40^3$ & 65 & 2.57e$-$1 & 39 & 1.66e$-$1 & 2.19 & 1.38 & 1.67 & 1.55 \\ \hline
$80^3$ & 127 & 3.85e+0 & 60 & 1.93e+0 & 2.42 & 1.48 & 2.12 & 2.00 \\ \hline
$160^3$ & 254 & 5.94e+1 & 98 & 2.39e+1 & 2.59 & 1.55 & 2.59 & 2.49 \\ \hline
$320^3$ & 503 & 9.21e+2 & 166 & 3.12e+2 & 2.72 & 1.59 & 3.03 & 2.95 \\ \hline
$640^3$ & 1015 & 1.59e+4 & 287 & 4.52e+3 & 2.80 & 1.62 & 3.54 & 3.52 \\ \hline
\end{tabular}
\end{center}
\label{table:result2}
\end{table}

\begin{table}[!h]
\caption{Reduction of the objective function by A$^2$ILU{\rm (0)} and its scalability.}
\begin{center} \footnotesize
\begin{tabular}{|c|c|c|c|} \hline
\multirow{2}{*}{$n$} & \multirow{2}{*}{ILU(0)} & \multirow{2}{*}{A$^2$ILU(0)} & Reduction \\
& {} & {} & ratio \\ \hline
$10^3$ & 4.16e+0 & 1.56e+0 & 0.38 \\ \hline
$20^3$ & 1.44e+1 & 3.77e+0 & 0.26 \\ \hline
$40^3$ & 4.36e+1 & 8.55e+0 & 0.20 \\ \hline
$80^3$ & 1.27e+2 & 1.86e+1 & 0.15 \\ \hline
$160^3$ & 3.66e+2 & 3.98e+1 & 0.11 \\ \hline
$320^3$ & 1.04e+3 & 8.38e+1 & 0.08 \\ \hline
$640^3$ & 2.96e+3 & 1.75e+2 & 0.06 \\ \hline
\end{tabular}
\end{center}
\label{table:result3}
\end{table}

From Table \ref{table:result2}, the speed-up ratio through auto-acceleration increases with matrix size, eventually exceeding 3.5 times the original performance.
The condition number of a linear system created by the digitization of second-order elliptic PDE, as we focused on in this article, is $O(h^{-2})$ regardless of the order of differentiation and the dimension of the space of interest.
The order of the number of iterations until convergence for CG with no preconditioning is $O(h^{-1})$.
This order estimation is still valid even if a simple preconditioning including ILU(0) is applied.
The order of condition number reduces to $O(h^{-1})$ and the number of iterations reduces to $O(h^{-0.5})$ when we apply MILU(0) preconditioning.
On the other hand, A$^2$ILU(0) improves the order of the number of iterations to $O(h^{-0.65})$.
Since our method is based on MILU methodology in constructing the objective function, it inherits the merits of MILU in the scalability with respect to the problem size.

Table \ref{table:result3} shows the objective function of ILU(0) and A$^2$ILU(0) and the reduction ratio defined as f(R) of A$^2$ILU(0) divided by that of ILU(0).
This table shows that as the matrix becomes larger, the objective function drastically decreases through auto-acceleration.
The decrease in the objective function indicates that the approximation accuracy of the preconditioning matrix is improved, as explained in section 3.

Hence, from Tables \ref{table:result2} and \ref{table:result3}, we conjecture that in ILU(0) preconditioning, the larger the coefficient matrix, the lower the approximation accuracy of the preconditioning matrix.
However, when auto-acceleration is applied, the acceleration parameters are optimized for each preconditioning matrix, regardless of the matrix size, thus improving the properties of the preconditioning matrix as much as possible.
As the matrix becomes larger, therefore, the effects of auto-acceleration also increase, improving the performance of A$^2$ILU(0) preconditioning relative to ILU(0) preconditioning.
For these reasons, we conclude that A$^2$ILU(0) preconditioning improves the scalability of ILU(0) preconditioning.

\subsection{Performance evaluation for general sparse matrices}
In this subsection we evaluate the performance of A$^2$ILU preconditioning for general linear systems without limiting the method for discretizing PDEs and validate its effectiveness.
In this case, the coefficient matrix for the system of equations is an irregular matrix in which nonzero entries are located irregularly.
We obtained our samples of coefficient matrices from the University of Florida sparse matrix collection under the conditions shown in Table \ref{table:selectionCriteria}.
To ensure the reliability of our evaluation, we carried out our validation with as many sample matrices as possible.

\begin{table}[!h]
\caption{Criteria for collecting sample matrices.}
\begin{center} \footnotesize
\begin{tabular}{|c|c|} \hline
Number of rows & $\le$ 600,000 \\ \hline
Number of nonzeros & $\le$ 20,000,000 \\ \hline
Pattern symmetry & Yes \\ \hline
Numerical symmetry & Yes \\ \hline
Shape & Square \\ \hline
Positive definite? & Either \\ \hline
2D/3D discretization? & Yes \\ \hline
Real or complex? & Real \\ \hline
Binary & No \\ \hline
\end{tabular}
\end{center}
\label{table:selectionCriteria}
\end{table}

For vector $b$ on the right-hand side of linear systems, we assigned values such that the solution vector $x$ was a vector each of whose entries was 1.
The initial value for each of the entries in the solution vector $x$ was 0.
In addition to the requirements shown in Table \ref{table:selectionCriteria}, we removed all matrices that contain zeros on the main diagonal and all matrices in which all the entries have such small absolute values that vector $b$ on the right-hand side becomes the zero vector; a total of 217 matrices were considered.
For the iterative method we used the CG method and for determination of convergence we used the criteria $\|r\|^2/\|b\|^2 \le$ 1.0e$-$8.
The maximum number of iterations was set as the size of the coefficient matrix.
For each problem we carried out diagonal scaling in advance to normalize the diagonal entries (making them 1) in the coefficient matrix.
Table \ref{table:calculationEnvironment2} summarizes the specifications of the computation environment.

\begin{table}[!h]
\caption{Calculation environment.}
\begin{center} \footnotesize
\begin{tabular}{|c|c|c|} \hline
System & Fujitsu PRIMERGY RX200 S3 \\ \hline
\multirow{2}{*}{CPU} & Intel Xeon 5160 \\
& 3.00 GHz / L2-Cache 4 MB \\ \hline
Memory & DDR2-667 8 GB \\ \hline
Compiler & Intel C++ Ver.11.1.075 \\ \hline
\end{tabular}
\end{center}
\label{table:calculationEnvironment2}
\end{table}

When solving systems of linear equations arising from an unstructured grid, ILU factorization may result in a diagonal matrix $D$ with some tiny entries, which significantly degrade the convergence \cite{shiftedILU}.
Shifted ILU preconditioning is an effective way to avoid this occurrence and thus is widely used, so we now validate the effectiveness of auto-acceleration on shifted ILU preconditioning.
To obtain the best performance of the original shifted ILU, we use a set of candidate values for the parameters.
We applied shifted ILU(0) preconditioning and shifted A$^2$ILU(0) preconditioning to every combination of the sample matrices and the candidate parameter values and analyzed the results from the following multiple aspects.

\subsubsection{Robustness of shifted A$^2$ILU(0)}
\begin{figure}[!h]
\begin{center}
\includegraphics[scale=\figscale,clip]{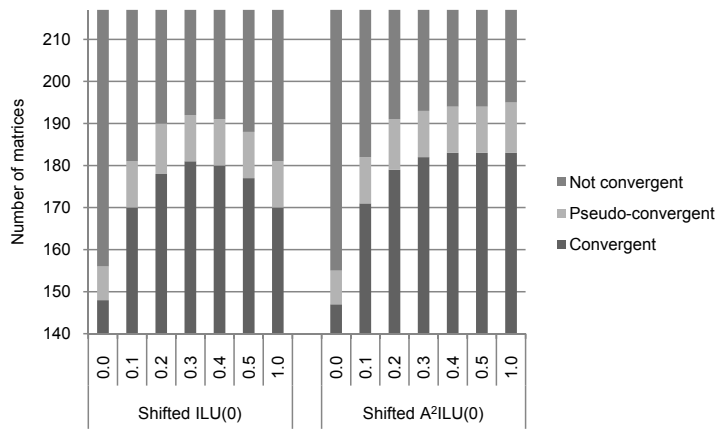}
\caption{Result of convergence determination.}
\label{fig:convergenceDetermination}
\end{center}
\end{figure}

In this subsection we examine how auto-acceleration affected the convergence.
Figure \ref{fig:convergenceDetermination} shows the total result of convergence determination for shifted ILU(0) preconditioning and shifted A$^2$ILU(0) preconditioning.
The left graph is the result of shifted ILU(0) and the right shifted A$^2$ILU(0).
The horizontal axis indicates the values of shift parameters, while the vertical axis indicates the number of matrices tallied.
The type of convergence is classified into the following three: convergent, pseudo-convergent, and not convergent.
Pseudo-convergent refers to a case where the norm of the true residual vector reaches its lower bound, while the norm of the recursive residual vector continues to decrease.
In this case, an approximate solution with a desired accuracy is not expected to be obtained even if the iteration proceeds.

As shown in the left side of the figure, the shifted ILU(0) made only 148 matrices converge where $\alpha = 0.0$.
The number of convergent cases increased up to 181 as $\alpha$ increased until $\alpha = 0.3$.
Hence, the results clearly show the effects of using shift parameters.
However, when $\alpha$ increased beyond 0.3, this tendency reversed itself, decreasing the number of convergent cases monotonically.
Thus, $\alpha=0.3$ is the optimum value for this collection of matrices.

Next, we look at the results of the shifted A$^2$ILU(0).
Overall, any reduction in the number of convergent cases compared to the shifted ILU(0) is not found for any value of $\alpha$.
When $\alpha = 0.4$ or above, there is a significant increase, avoiding the reduction in the number of convergent cases with an increase of the value of $\alpha$.
Hence, the use of auto-acceleration keeps or enhances the robustness of shifted ILU(0) preconditioning against the variety of coefficient matrices.

\subsubsection{Improved convergence rate by auto-acceleration}
In this subsection we further analyze the effects of auto-acceleration, not only on the types of convergence but also the speed of convergence.
In the discussion below, we evaluate the computing time by using the number of iterations because we found that the computing cost of auto-acceleration can be ignorable for the total cost (no more than 1 $\%$ for $\alpha = 0.2, 0.5$).

\begin{figure}[!h]
\begin{center}
\includegraphics[scale=\figscale,clip]{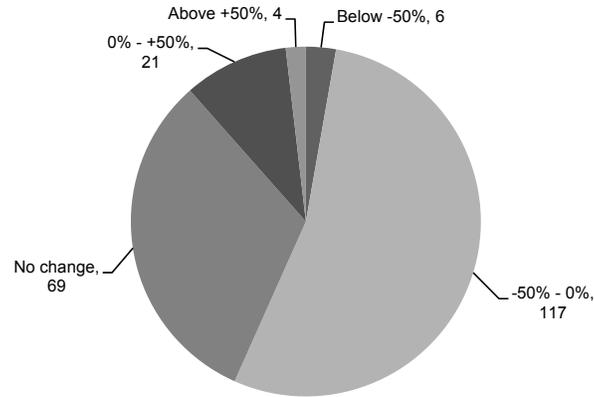}
\caption{Increase ratio of the iterations through auto-acceleration (shift parameter $\alpha = 0.2$).}
\label{fig:convergenceRate0.2}
\end{center}
\end{figure}

\begin{figure}[!h]
\begin{center}
\includegraphics[scale=\figscale,clip]{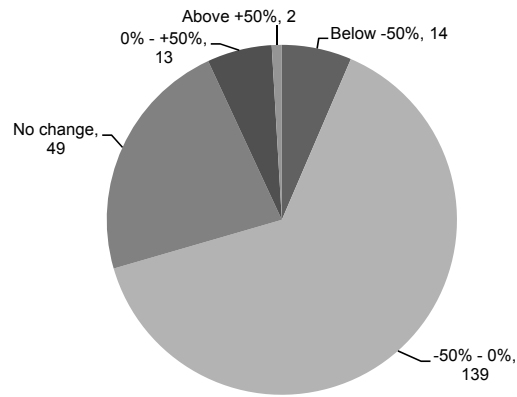}
\caption{Increase ratio of the iterations through auto-acceleration (shift parameter $\alpha = 0.5$).}
\label{fig:convergenceRate0.5}
\end{center}
\end{figure}

We have divided the increase ratio of the number of iterations into several classes and show the number of matrices in each of these classes in Figures \ref{fig:convergenceRate0.2} and \ref{fig:convergenceRate0.5}.
Here, the increase ratio is defined as $(N_A - N_I) / N_I$, where $N_A$ and $N_I$ are the number of iterations obtained in shifted A$^2$ILU(0) and shifted ILU(0), respectively.
Figure \ref{fig:convergenceRate0.2} shows the results when the shift parameter is $\alpha = 0.2$, while Figure \ref{fig:convergenceRate0.5} shows the results at $\alpha = 0.5$.
We assigned the ratio of below $-$50\% to a problem that converges by shifted A$^2$ILU(0) but not by shifted ILU(0).
The inverse case was indicated by above $+$50\%.
There is no change if both methods obtained convergence by the same number of iterations or if neither of them produced convergence.

Figure \ref{fig:convergenceRate0.2} shows that when shift parameter $\alpha$ is 0.2, the number of solutions whose convergence rate improves by auto-acceleration is 123 (below $-$50\% and $-$50\% -- 0\%), which is 56.7\% overall.
Meanwhile, the convergence rate remained the same for 69 solutions (no change), or 31.8\% of the total number.
The number of solutions whose convergence rate got worse was only 25 (0\% -- $+$50\% and above $+$50\%), which is 11.5\% of the total number.
The application of the auto-acceleration resulted in many more merits than demerits in convergence rate as well as robustness.

Further, Figure \ref{fig:convergenceRate0.5} shows that at shift parameter $\alpha = 0.5$, the percentage in which the convergence rate is improved by the application of auto-acceleration increases to 70.5\%, while the convergence rate did not change in 22.6\% and got worse in 6.9\% (both of these numbers dropped here).
Hence, we conclude that the auto-acceleration more clearly shows its effectiveness with a larger shift parameter value $\alpha$.

\subsubsection{Effect of shifted A$^2$ILU(0) for the shift parameter}
\begin{figure}[!h]
\begin{center}
\includegraphics[scale=\figscale,clip]{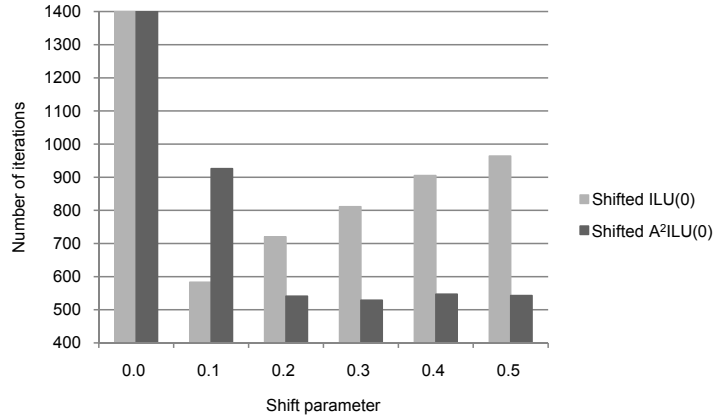}
\caption{Effects of auto-acceleration in Rothberg/cfd2.}
\label{fig:cfd2}
\end{center}
\end{figure}

In subsections 4.2.1 and 4.2.2, we demonstrated the effectiveness of auto-acceleration on the robustness over the variety of general sparse matrices from the standpoint of performance.
In this subsection, we further show the effectiveness of auto-acceleration on the robustness over the shift parameter values from the standpoint of applicability.
Specifically, we address the selection of the shift parameter, a major practical challenge in shifted ILU(0), by studying what effects the auto-acceleration has on parameter selection and verifying the practicality of shifted A$^2$ILU(0).
To help us understand this, we show in Figure \ref{fig:cfd2} the calculation results for the sample matrix Rothberg/cfd2.
The horizontal axis represents the shift parameter while the vertical axis represents the number of iterations.

There are two main effects of using a shift parameter in shifted ILU preconditioning.
One is the positive effect of helping the convergence of solutions that would otherwise not converge.
Figure \ref{fig:cfd2} shows this effect when $\alpha = 0.1$.
The other is the negative effect of gradually reducing the convergence rate.
This is shown by the increase of the number of iterations when $\alpha$ exceeds 0.1.
Because of these two opposing effects, the responsiveness of shifted ILU(0) to the shift parameter value is quite sensitive.
In contrast, with shifted A$^2$ILU(0) preconditioning, only the latter effect, namely, the negative effect of reduced convergence rate when $\alpha$ exceeds 0.1, is kept low.
So the robustness of shifted A$^2$ILU(0) against the shift parameter value is improved.
We examine this point below.

The positive effect of the shift parameter exists for the following reason.
In ILU preconditioning, the convergence is significantly degraded if, during ILU factorization, a tiny entry appears in diagonal matrix $D$.
Shifted ILU avoids this problem by enlarging the diagonal entries of the matrix being ILU factorized.
As a result, linear systems whose solutions are not supposed to converge can have convergent solutions.
The drastic increase, seen in Figure \ref{fig:convergenceDetermination}, in the number of convergent cases for $\alpha = 0.0$ -- $0.2$ is thought to be caused by this effect.
On the other hand, the auto-acceleration process does not change the sign of each entry of diagonal matrix $D$ (even though it scales the matrix) and so does not prevent this effect.
Therefore, shifted A$^2$ILU(0) has proved as effective as shifted ILU(0), a fact revealed by Figures \ref{fig:cfd2} and \ref{fig:convergenceDetermination}.

Next, the negative effect of the shift parameter exists for the following reason.
In shifted ILU, the larger the parameter is, the farther the matrix being ILU factorized becomes from coefficient matrix $A$.
This reduces the accuracy of the preconditioning matrix, causing the convergence rate to decrease.
The decrease in the number of convergent cases when $\alpha$ is 0.3 or greater, shown in Figure \ref{fig:convergenceDetermination}, is thought to be caused by this effect.
On the other hand, the auto-acceleration process improves this accuracy of the preconditioning matrix by the acceleration parameters, reducing this negative effect.
In particular, because the acceleration parameters are automatically optimized in accordance with the remainder matrix, this effect becomes more pronounced as the shift parameter increases.
This fact is supported by the results of Figures \ref{fig:convergenceRate0.2} and \ref{fig:convergenceRate0.5} that the larger the value of $\alpha$, the higher the percentage where the convergence rate improved.
Consequently, as Figures \ref{fig:cfd2} and \ref{fig:convergenceDetermination} show, shifted A$^2$ILU(0) maintains high performance even when the shift parameter takes larger values than the optimum one.
For these reasons, we conclude that the auto-acceleration cancels only the negative effect of the shift parameter while smoothing out the responsiveness to shift parameter increases.

With this stated, we consider the practicality of shifted ILU(0) and shifted A$^2$ILU(0).
Shifted ILU is often used to avoid the breakdown caused by tiny diagonal elements.
As described above, a user needs to find the minimum value of $\alpha$ that avoids the breakdown.
To do this, shifted ILU(0) requests a user to perform a line search along the value of $\alpha$, in which one should check whether the breakdown is occurs at each point on the line.
This brute force optimization consumes obviously extremely huge computational cost.
Shifted A$^2$ILU(0) is expected to request no such brute force optimization of $\alpha$ because A$^2$ inhibits the negative effects.
Therefore, shifted A$^2$ILU(0) tremendously improves the practicality of shifted ILU(0).

\section{Conclusions}
In this paper we have proposed A$^2$ILU preconditioning, which improves performance without losing the practicality of ILU preconditioning.
A$^2$ILU preconditioning is a process in which new acceleration parameters are incorporated in ILU preconditioning and these parameters are automatically optimized.
Previous ILU-based methods all have practicality issues because their own parameters must be set up by users.
In contrast, A$^2$ILU(0) preconditioning is highly practical because, like ILU(0) preconditioning, it is a parameter-free method for users.

We verified the following merits of A$^2$ILU(0) preconditioning by using systems of linear equations arising from physical simulations based on rectangular grids:
\begin{itemize}
\item For five sample problems with the coefficient matrices of hundred-thousands dimension, A$^2$ILU(0) preconditioning is 1.65 times as fast as ILU(0) preconditioning on average.
Even compared with other ILU-based methods in which their original parameters are optimized manually, this speed is still higher.
We concluded that A$^2$ILU(0) is superior to the previous ILU-based methods with respect to both practicality and performance.
\item Its scalability relative to the size of the coefficient matrix is superior to ILU(0) preconditioning.
It is concluded that the number of iterations is ${\cal O}(h^{-0.65})$ which is less than ${\cal O}(h^{-1})$ for ILU(0), where $h$ denotes the mesh size.
Once the matrix size exceeds 260 million, its speed exceeds 3.5 times that of ILU(0) preconditioning.
\end{itemize}

Furthermore, the proposed auto-acceleration was applied to previous major ILU-based methods with the following confirmed merits:
\begin{itemize}
\item For shifted ILU(0), modified ILU(0), and fill-in level ILU, regardless of the value of the parameter unique to these methods, the performance improves.
\item For Crout ILU, the performance improves even in the case when a denser preconditioned matrix is generated and therefore the original Crout ILU shows good performance.
\item For shifted ILU(0) and modified ILU(0), because the performance is stable with respect to any change in the unique parameter, the burden of setting up the parameter is reduced so that its practicality is improved.
\end{itemize}

For general sparse matrices, we have shown the effectiveness of shifted A$^2$ILU(0).
We evaluated the performance of shifted A$^2$ILU(0) preconditioning using over 200 general sparse matrices obtained from the University of Florida sparse matrix collection.
The results confirmed the following merits:
\begin{itemize}
\item There is no reduction in the number of convergent cases compared with shifted ILU(0) preconditioning over all the shift parameters examined in this paper.
In addition, when the shift parameter is beyond 0.3, auto-acceleration increased the number of convergent cases of shifted ILU(0).
\item Many more cases improve the convergence rate, rather than worsen, compared with shifted ILU(0) preconditioning.
\item Even if the value of the shift parameter is raised to ensure convergence, the convergence rate does not drop significantly, unlike shifted ILU(0) preconditioning.
Hence, this method is able to maintain both safety and performance and is thus more practical.
\end{itemize}

\section*{Acknowledgement}
The authors would like to thank anonymous reviewers for their insightful comments and encouragement.

\end{document}